\documentclass[12pt]{amsart}
\usepackage{amsmath}
\usepackage{amssymb}
\usepackage{epsfig}
\newtheorem{theorem}{Theorem}
\newtheorem{lemma}{Lemma}
\newtheorem{corollary}{Corollary}
\DeclareMathOperator{\Res}{Res}

\DeclareMathOperator{\Jac}{Jac}
\DeclareMathOperator{\Prym}{Prym}
\DeclareMathOperator{\Nm}{Nm}
\newcommand{\Z}{\mathbb{Z}}
\newcommand{\C}{\mathbb{C}}

\usepackage{graphics}

\begin{document}

\title{Dimensions of Prym Varieties}
\author{Amy E. Ksir\\
        Mathematics Department\\
        State University of New York at Stony Brook \\
        Stony Brook, NY, 11794}
\date{July 26, 2000}
\email{ksir@math.sunysb.edu}
\begin{abstract}
Given a tame Galois branched cover of curves $\pi: X \to Y$ with any finite
Galois group $G$ whose representations are rational, we compute the dimension 
of the (generalized) Prym variety $\Prym_{\rho}(X)$ corresponding to any 
irreducible representation $\rho$ of $G$.  This formula can be applied to the 
study of algebraic integrable systems using Lax pairs, in particular systems
associated with Seiberg-Witten theory.  However, the formula is much 
more general and its computation and proof are entirely algebraic.
\end{abstract}
\maketitle

\section{Introduction}

The most familiar Prym variety arises from a (possibly 
branched) double cover $\pi: X \to Y$ of curves.  In this 
situation, there is a surjective norm map $\Nm: \Jac(X) \to 
\Jac(Y)$, and the Prym (another abelian variety) is a connected 
component of its kernel.  Another way to think of this is that the 
involution $\sigma$ of the double cover induces an action of $\Z/2\Z$ on 
the vector space $H^{0}(X, \omega_{X})$, which can then be 
decomposed as a representation of $\Z/2\Z$.  The Jacobian of the base 
curve Y and the Prym correspond to the trivial and sign 
representations, respectively.  The Prym variety can 
be defined as the component containing the identity of 
$(\Jac(X) \otimes_{\Z} \varepsilon)^{\sigma}$, where $\varepsilon$ denotes 
the sign representation of $\Z/2\Z$.

The generalization of this construction that we will study in this 
paper is as follows. Let $G$ be a finite group, and $\pi: X \to Y$ be a 
tame Galois branched cover, with Galois group $G$, of smooth projective 
curves over an algebraically closed field.  The action of $G$ on $X$ induces 
an action on the vector space of differentials $H^0(X,\omega_{X})$,
and on the Jacobian  $\Jac(X)$.  For any representation $\rho$ of 
$G$, we define $\Prym_{\rho}(X)$ to be the connected component 
containing the identity of $(\Jac(X) \otimes_{\Z} \rho^{*})^{G}$.
The vector space $H^0(X,\omega_{X})$ will decompose as a $\Z[G]$-module into a
direct sum of isotypic pieces
\begin{equation}
H^0(X,\omega_{X}) = \bigoplus_{j=1}^{N} \rho_j \otimes V_j
\end{equation}
where $\rho_1, \ldots , \rho_N$ are the irreducible representations of $G$.
If $G$ is such that all of its representations are rational, then the
Jacobian will also decompose, up to isogeny, into a direct sum of 
Pryms [D2]:
\begin{equation}
\Jac(X) \sim \bigoplus_{j=1}^{N} \rho_j \otimes \Prym_{\rho_j}(X).
\end{equation}
In particular, if $G$ is the Weyl group of a semisimple Lie algebra, then it 
will satisfy this property.

The goal of this paper is to compute the dimension of such a Prym 
variety.  This formula is given in section 2, with a proof that uses 
only the Riemann-Hurwitz theorem and some character theory.  Special
cases of this formula relevant to integrable systems have appeared 
previously [A, Me, S, MS].

One motivation for this work comes from the study of algebraically integrable 
systems.  An algebraically integrable system is a Hamiltonian system 
of ordinary differential equations, where 
the phase space is an algebraic variety with an algebraic (holomorphic, 
over $\C$) 
symplectic structure.  The complete integrability of the system means 
that there are $n$ commuting Hamiltonian functions on the 
$2n$-dimensional phase 
space.  For an algebraically integrable system, these functions should 
be algebraic, in which case they define a morphism to an 
$n$-dimensional space of states for the system.  The flow of the 
system will be linearized on the fibers of this morphism, which, if 
they are compact, will be $n$-dimensional abelian varieties.

Many such systems can be solved by expressing the system as a Lax
pair depending on a parameter $z$.  The equations can be written in the 
form $\frac{d}{dt}A = [A,B]$, where $A$ and $B$ are elements of a Lie
algebra $\mathfrak{g}$, and depend both on time $t$ and on a 
parameter $z$, which is thought of as a coordinate on a curve $Y$.
In this case, the flow of the system is linearized on a subtorus of
the Jacobian of a Galois cover of $Y$.  If it can be shown that this
subtorus is isogenous to a Prym of the correct dimension, then the
system is completely integrable.

In section 3, we will 
briefly discuss two examples of such systems, the periodic Toda lattice 
and Hitchin systems.  Both of these are important in Seiberg-Witten 
theory, providing solutions to $\mathcal{N}=2$ supersymmetric 
Yang-Mills gauge theory in four dimensions.

This work appeared as part of a Ph.D. thesis at the University of
Pennsylvania.  The author would like to thank her thesis advisor, Ron Donagi, 
for suggesting this project and for many helpful discussions.  Thanks
are also due to David Harbater, Eyal Markman, and Leon Takhtajan.

\section{Dimensions}

We can start by using the Riemann-Hurwitz formula to find the genus $g_X$ of 
$X$, which will be the dimension of the whole space $H^0(X,\omega_{X})$ and of
$\Jac (X)$.  Since $\pi: X \to Y$ is a cover of degree $|G|$,
we get
\begin{equation}
g_X = 1 + |G|(g-1) + \frac{\deg R}{2}
\end{equation}
where $g$ is the genus of the base curve $Y$ and $R$ is the ramification 
divisor.  

The first isotypic piece we can find the dimension of is $V_1$, corresponding
to the trivial representation.  The subspace where $G$ acts trivially is
the subspace of differentials which are pullbacks by $\pi$
of differentials on $Y$.   This tells us that $\dim V_1 = \dim 
H^0(Y,\omega_{Y}) = g$.

In the case of classical Pryms, where $G = \mathbb{Z}/2$, there is only one
other isotypic piece, $V_{\varepsilon}$ corresponding to the sign 
representation $\varepsilon$.  Thus we have
\begin{equation}
\dim V_{\varepsilon} = g_X - g = g-1 + \frac{\deg R}{2}.
\end{equation}

For larger groups $G$, there are more isotypic pieces, but we also have
more information:  we can look at intermediate
curves, i.e. quotients of $X$ by subgroups $H$ of $G$.  Differentials
on $X/H$ pull back to differentials on $X$ where $H$ acts trivially.  Thus 
\begin{equation}
H^0(X/H,\omega_{X/H}) = \bigoplus_{j=1}^{N} (\rho_j)^{H} \otimes V_j.
\end{equation}
The map $\pi_H: X/H \to Y$ will be a cover of degree $\frac{|G|}{|H|}$, so
Riemann-Hurwitz gives us the following formula for the genus $g_H$ of 
$X/H$, which is the dimension of $H^0(X/H,\omega_{X/H})$:
\begin{equation}
g_H = 1 + \frac{|G|}{|H|}(g-1) + \frac{\deg R_H}{2}.
\end{equation}
where again $R_H$ is the ramification divisor.

We can further analyze the ramification divisor, by classifying the branch
points according to their inertial groups.   Since 
$\pi: X \to Y$ is a Galois cover of curves over $\mathbb{C}$, all of the 
inertial groups must be cyclic.  

\begin{lemma}
Let $G$ be a finite group all of whose characters are defined over $\mathbb{Q}$.
If two elements $x,y \in G$ generate conjugate cyclic subgroups, then they are
conjugate.
\end{lemma}

Proof (adapted from [BZ]):  We want to show that for any character $\chi$ of 
$G$,  $\chi(x) = \chi(y)$.  Then the properties of characters will tell us 
that $x$ and $y$ must be in the same conjugacy class.  

We may assume that $x$ and $y$ generate the same subgroup, $H$.  Then 
$y = x^k$ for some integer $k$ relatively prime to $|H|$.  Let $\chi$
be a character of $G$, and $\rho: G \to GL(n,\mathbb{C})$ 
a representation with character $\chi$.  Then $\rho(x)$ will be a matrix
with eigenvalues $\lambda_1, \ldots, \lambda_n$, and $\rho(y)$ will
have eigenvalues $\lambda_1^k, \ldots, \lambda_n^k$.  Since $x^{|H|} = 1$, 
we have $\lambda_1^{|H|} = \ldots = \lambda_n^{|H|} = 1$.  Let $\xi$ be
a primitive $|H|$th root of unity.  Then we can write $\lambda_1 = 
\xi^{\nu_1}, \ldots, \lambda_n = \xi^{\nu_n}$ for some integers $\nu_i$.
Now $\chi(x)$ = Trace($\rho(x)$) = $\lambda_1 +  \ldots + \lambda_n$,
and $\chi(y)$ = $\chi(x^k)$ = $\lambda_1^k +  \ldots + \lambda_n^k$.
Thus $\chi(y)$ will be the image of $\chi(x)$ under the element of
Gal($\mathbb{Q}(\xi)/\mathbb{Q}$) which sends $\xi \mapsto \xi^k$.
Since the values of $\chi$ are rational, this element will act trivially,
so $\chi(y) = \chi(x)$. \hfill $\square$ \smallskip

 From now on, we will suppose that $G$ is such that all of its characters are
rational.  (This will be true, for instance, if $G$ is a Weyl group).  Pick 
representative
elements $h_1 \ldots h_N$ for each conjugacy class in $G$,
and let $H_1 \ldots H_N$ be the cyclic groups that each of them generates.
By Lemma 1, this will be the whole set (up to conjugacy) of cyclic 
subgroups of $G$.  We can partially order this set of cyclic subgroups
by their size, so that $H_1$ is the trivial subgroup. Now we can classify 
the branch points:  let $R_k, k=2 \ldots N$ be the degree of the branch 
locus with inertial group conjugate to $H_k$ (ignoring the trivial group).  
Over each point of the branch
locus where the inertial group is conjugate to $H_k$, there will be $|G|/|H_k|$
points in the fiber. Thus the degree of the ramification 
divisor $R$ of  $\pi: X \to Y$ will be
\begin{equation}
\deg R = \sum_{k=1}^{N} (|G| - \frac{|G|}{|H_k|}) R_k
\end{equation}
For each quotient curve $X/H$, each point in the fiber of $\pi_H: X/H \to Y$ 
over a point with inertial group $H_k$ will correspond to a double coset
$H_k \backslash G / H$.  Thus the degree of the ramification divisor $R_H$
will be
\begin{equation}
\deg R_H = \sum_{k=1}^{N} (\frac{|G|}{|H|} - \#(H_k \backslash G / H)) R_k.
\end{equation}
Combining these formulas with the earlier Riemann-Hurwitz computations, we get:
\begin{equation}
g_X = 1 + |G|(g-1) + \sum_k(|G| - \frac{|G|}{|H|} ) \frac{R_k}{2}
\end{equation}
\begin{equation}
g_H = 1 + \frac{|G|}{|H|}(g-1) + \sum_k(\frac{|G|}{|H|} - 
\#(H_k \backslash G / H) ) \frac{R_k}{2}
\end{equation}

Since the genera $g_H$ are exactly the dimensions $\dim 
H^0(X/H,\omega_{X/H})$,
we also have 
\begin{equation}
g_H = \sum_{j=1}^{N} \dim \rho_j^{H} \dim V_j.
\end{equation}

For each subgroup $H$, this is a linear equation for 
the unknown dimensions $\dim V_j$ in terms of the genus $g_{H}$.
Thus by taking quotients by the set of all cyclic subgroups $H_1 \ldots H_N$, 
we get a system of $N$ equations.  We wish to 
invert the matrix $\dim\rho_j^{H_i}$ and find the $N$ unknowns $\dim V_j$.
\begin{lemma}
The matrix $\dim\rho_j^{H_i}$ is invertible.
\end{lemma}

Proof:  We show that the rows of the matrix are linearly independent, using 
the fact that rows of the character table are linearly independent. 
First, note that $\dim \rho_j^{H_i}$, the dimension of the subspace of $\rho_j$
invariant under $H_i$, is equal to the inner product of characters
$\langle \Res^G_{H_i} \rho_j, \mathbf{1} \rangle$, 
which we can read off from the character table of $G$ as
\begin{equation}
\dim \rho_j^{H_i} = \frac{1}{|H_i|} \sum_{a_i \in H_i} \chi_{\rho_j}(a_i).
\end{equation}

Compare this matrix to the matrix of the character table $\chi_{\rho_j}(a_i)$.
From (12) we see that each row is a sum of multiples of rows of the 
character table.  Since each element of a subgroup has order less than or
equal to the order of the subgroup, the rows of the
character table being added to get row $i$ appear at or below row $i$ in
the character table.  Thus if we write the matrix $\dim\rho_j^{H_i}$
in terms of the basis of the character table, we will get a lower triangular
matrix with non-zero entries on the diagonal.  By row reduction, we see
that the linear independence of the rows of $\dim\rho_j^{H_i}$ is equivalent
to the linear independence of the rows of the character table. \hfill $\square$

\begin{theorem}
For each nontrivial irreducible representation $\rho_j$ of $G$, 
$V_j$ has dimension
\begin{equation}
(\dim \rho_j) (g-1) + \sum_{k=1}^{N} \Bigl((\dim \rho_j) - (\dim \rho_j^{H_k})
        \Bigr) \frac{R_{H_k}}{2}
\end{equation}
\end{theorem}

Proof:  Since the matrix $\dim\rho_j^{H_i}$ is invertible, there is a 
unique solution to the system of equations (11), so we only need to show
that this is a solution.  Namely, given this formula for $\dim V_j$,
and combining (10) and (11),
we wish to show that for each cyclic subgroup $H_i$,
\begin{equation}
\sum_{j=1}^{N} \dim \rho_j^{H_i} \dim V_j = 1 + \frac{|G|}{|H_i|}(g-1) 
+ \sum_k(\frac{|G|}{|H_i|} - \#(H_k \backslash G / H_i) ) \frac{R_k}{2}.
\end{equation}

Note that on the left side we are summing over all representations, not just
the nontrivial ones, so our notation will be simpler if we write 
$\dim V_1 = g$ in a similar form to (11).  For the trivial representation 
$\rho_1$, $(\dim \rho_1) - (\dim \rho_1^{H_k}) = 0$ (since $\rho_1$ is fixed by 
any subgroup $H_k$), so 
\begin{equation}
\dim V_1 = 1 + (\dim \rho_1) (g-1) + \sum_{k=1}^{N} \Bigl((\dim \rho_1) - 
	(\dim \rho_1^{H_k}) \Bigr) \frac{R_{H_k}}{2}.
\end{equation}

The sum on the left hand side of (14) will be
\begin{equation}
1 + \sum_{j=1}^{N} \dim \rho_j^{H_i} \Bigl((\dim \rho_j) (g-1) 
+ \sum_{k=1}^{N} \bigl((\dim \rho_j) - (\dim \rho_j^{H_k}) \bigr) 
\frac{R_{H_k}}{2} \Bigr). 
\end{equation}

Let us look at the $(g-1)$ term and the $R_{H_k}$ terms separately.
For the $(g-1)$ coefficient, we can write both $\dim \rho_j^{H_i}$  and 
$\dim \rho_j$ in terms of characters of $G$ (as in (12)) and exchange the 
order of summation to get

\begin{equation}
\sum_{j=1}^{N} \dim \rho_j^{H_i} \dim \rho_j =
 \frac{1}{|H_i|} \sum_{a_i \in H_i} \sum_{j=1}^N \chi_{\rho_j}(a_i) 
	\chi_{\rho_j}(e) 
\end{equation}
where $e$ is the identity element of $G$.  The inner sum amounts to taking
the inner product of two columns of the character table of $G$.  The
orthogonality of characters tells us that this inner product will be zero
unless the two columns are the same, in this case if $a_i = e$.
Thus the sum over elements in $H_i$ disappears, and we get the sum of the 
squares of the dimensions of the characters:
\begin{equation}
 \frac{1}{|H_i|} \sum_{j=1}^N \chi_{\rho_j}(e)^2 = \frac{|G|}{|H_i|}.
\end{equation}
which is what we want.

The $R_{H_k}$ term looks like
\begin{equation}
\sum_{j=1}^{N} \dim \rho_j^{H_i} \sum_{k=1}^{N} \bigl((\dim \rho_j) - 
	(\dim \rho_j^{H_k}) \bigr) \frac{R_{H_k}}{2}.
\end{equation}

We can distribute and rearrange the sums to get:
\begin{equation}
\sum_{k=1}^{N} \Bigl(\sum_{j=1}^{N} \dim \rho_j^{H_i} \dim \rho_j -
	\sum_{j=1}^{N} \dim \rho_j^{H_i} \dim \rho_j^{H_k} 
	\Bigr)\frac{R_{H_k}}{2}.
\end{equation}

As in (17) and (18), the first term becomes $\frac{|G|}{|H_i|}$.  The 
second term is also the inner product of columns of the character table:
\begin{equation}
\sum_{j=1}^{N} \dim \rho_j^{H_i} \dim \rho_j^{H_k} =
	\frac{1}{|H_i|} \frac{1}{|H_k|} \sum_{a_i \in H_i} \sum_{a_k \in H_k}
	\sum_{j=1}^{N} \chi_{\rho_j}(a_i) \chi_{\rho_j}(a_k).
\end{equation}
This will be 
zero unless $a_i$ and $a_k$ are conjugate, in which case $\chi_{\rho_j}(a_i)  
= \chi_{\rho_j}(a_k)$ and character theory tells us (see for example [FH],
p. 18) that
\begin{equation}
\sum_{j=1}^N \chi_{\rho_j}(a_i)^2 = \frac{|G|}{c(a_i)},
\end{equation}
where $c(a_i)$ is the number of elements in the conjugacy class of $a_i$.  
Now the second term has become
\begin{equation}
\frac{|G|}{|H_i||H_k|} \sum_{\{a_i, a_k\}} \frac{1}{c(a_i)}
\end{equation}
where the sum is taken over pairs of elements $a_i \in H_i, a_k \in H_k$
such that $a_i$ and $a_k$ are conjugate.  This is exactly the number of
double cosets  $\#(H_k \backslash G / H_i)$.

Adding up all of the terms, the sum on the left hand side becomes
\begin{equation}
1 + \frac{|G|}{|H_i|}(g-1) + (\frac{|G|}{|H_i|} - \#(H_k \backslash G / H_i))
	\frac{R_{H_k}}{2}
\end{equation}
which is exactly the right hand side. \hfill $\square$ \smallskip

\begin{corollary}
For each nontrivial irreducible representation $\rho_j$ of $G$,
$\Prym_{\rho_j}(X)$ has dimension
\begin{equation}
(\dim \rho_j) (g-1) + \sum_{k=1}^{N} \Bigl((\dim \rho_j) - (\dim \rho_j^{H_k})
        \Bigr) \frac{R_{H_k}}{2}.
\end{equation}
\hfill $\square$
\end{corollary}

\section{Integrable Systems.}

\textbf{Periodic Toda lattice.}

The periodic Toda system is a Hamiltonian system of differential 
equations with Hamiltonian 
\begin{equation*}
    H(p,q) = \frac{|p|^{2}}{2} + \sum_{\alpha} e^{\alpha(q)}
\end{equation*}
where $p$ and $q$ are elements of the Cartan subalgebra $\mathfrak{t}$ 
of a semisimple Lie algebra $\mathfrak{g}$, and the sum is over the 
simple roots of $\mathfrak{g}$ plus the highest root.  This system 
can be expressed in Lax form [AvM] $\frac{d}{dt}A = [A,B]$, where $A$ 
and $B$ are elements of the loop algebra $\mathfrak{g}^{(1)}$.  and 
can be thought of as elements of $\mathfrak{g}$ which depend on a 
parameter $z \in \mathbb{P}^{1}$.  For $\mathfrak{sl}(n)$, $A$ is of 
the form
\begin{equation*}
\begin{pmatrix}
    y_{1} & 1     &   & x_{0}z \\
    x_{1} & y_{2}  & \ddots & \\
          & \ddots &  \ddots & 1 \\
    z     &        & x_{n-1} & y_{n}
\end{pmatrix}
\end{equation*}

For any representation $\varrho$ of 
$\mathfrak{g}$, the spectral curve $S_{\varrho}$ defined by the 
equation $\det (\varrho(A(z) - \lambda I) = 0$ is independent of 
time (i.e. is a conserved quantity of the system).  The spectral curve is a 
finite cover of $Y$ which for generic $z$ parametrizes the eigenvalues of
$\varrho(A(z))$.  While the eigenvalues are conserved by the system, 
the eigenvectors are not.  The eigenvectors of $\varrho(A)$ determine 
a line bundle on the spectral cover, so an element of 
$\Jac(S_{\varrho})$.  The flow of the system is linearized on this Jacobean.  
Since the original system of equations didn't depend on a choice of 
representation $\varrho$, the flow is actually linearized on an abelian 
variety which is a subvariety of $\Jac(S_{\varrho})$ for every 
$\varrho$.  

In fact, instead of considering each spectral cover we can look at 
the cameral cover $X \to \mathbb{P}^{1}$.  This is constructed as a pullback to 
$\mathbb{P}^{1}$ of the cover
$\mathfrak{t} \to \mathfrak{t}/G$, where $G$ is the Weyl group of 
$\mathfrak{g}$.  This cover is pulled back by the rational map
$\mathbb{P}^{1} \dashrightarrow \mathfrak{t}/G$ defined by the class of $A(z)$ 
under the adjoint action of the corresponding Lie group.  (For $A(z)$ 
a regular semisimple element of $\mathfrak{sl}(n)$, this map sends 
$z$ to the unordered set of eigenvalues of $A(z)$.)  Thus the cameral 
cover is a finite Galois cover of $\mathbb{P}^{1}$ whose Galois group $G$ 
is the Weyl group of $\mathfrak{g}$.  The flow of the Toda system is 
linearized on the Prym of this cover corresponding to the representation 
of $G$ on $\mathfrak{t}^{*}$.  This is an $r$-dimensional representation, 
where $r$ is the rank, so the dimension of this Prym is
\begin{equation*}
    r(-1) + \sum_{k=1}^{N} \Bigl(r - (\dim \mathfrak{t}^{H_k})
        \Bigr) \frac{R_{H_k}}{2}.
\end{equation*}

The ramification of this cover has been analyzed in [D1] and [MS].  
There are $2r$ branch points where the inertial group $H$ is $\Z/2\Z$ 
generated by one reflection, so for each of these $\dim 
\mathfrak{t}^{H}$ is $r-1$.  There are also two points ($z=0$ and 
$\infty$) where the inertial group $H$ is generated by the Coxeter 
element, the product of the reflections corresponding to the simple 
roots.  This element of $G$ doesn't fix any element of 
$\mathfrak{t}$, so for these two points $\dim \mathfrak{t}^{H} = 0$.
Thus the dimension of the Prym is
\begin{eqnarray*}
    -r + (r-(r-1)) \frac{2r}{2} + (r-0) \frac{2}{2} \\
    = r.
\end{eqnarray*}
Since the original system of equations had a $2r$-dimensional phase 
space, this is the answer that we want.

\medskip

\textbf{Hitchin systems.}  Hitchin showed [H] that the cotangent 
bundle to the moduli space of semistable
vector bundles on a curve $Y$ has the structure of an algebraically 
completely integrable system.  His proof, later extended to 
principal $\mathcal{G}$ bundles with any reductive Lie group $\mathcal{G}$ [F,S],
uses the fact that this moduli space is equivalent (by deformation
theory) to the space of \emph{Higgs pairs}, pairs $(P,\phi)$ of a
principal bundle and an endomorphism $\phi \in H^{0}(Y, ad(P) \otimes 
\omega_{Y})$.   
As in the case of the Toda system, the key construction is of a 
cameral cover of $Y$.  The eigenvalues of $\phi$, which are sections 
of the line bundle $\omega_{Y}$, determine a spectral cover of $Y$ in the 
total space bundle.  The eigenvectors determine a line bundle on this 
spectral cover.  The Hitchin map sends a Higgs pair $(P,\phi)$ to 
the set of coefficients of the characteristic polynomial.  Each coefficient is a 
section of a power of $\omega_{Y}$, so the image of the Hitchin map is
$B := \bigoplus_{i=1}^{r}H^{0}(Y,\omega_{Y}^{\otimes d_i})$, where 
the $d_{i}$ are the degrees of the basic invariant polynomials of the Lie 
algebra $\mathfrak{g}$.

Again, we can consider instead the cameral cover $X_{b} \to Y$, 
which is obtained as a pullback to $Y$ vi $\phi$ of $\mathfrak{t} \otimes 
\omega_{Y} \to \mathfrak{t} \otimes \omega_{Y}/G$.
The generic fiber of the Hitchin map is isogenous to 
$\Prym_{\mathbb{t}}(X)$, which has dimension
\begin{equation*}
    r(g-1) + \sum_{k=1}^{N} \Bigl(r - (\dim \mathfrak{t}^{H_k})
        \Bigr) \frac{R_{H_k}}{2}
\end{equation*}
By looking at the generic fiber,  
we can restrict our attention to cameral covers where the only 
ramification is of order two, with inertial group $H$ generated by one 
reflection.  The last piece of information we need to compute the 
dimension is the degree of the branch divisor of $X \to Y$.
 
The cover $\mathfrak{t} \otimes 
\omega_{Y} \to \mathfrak{t} \otimes \omega_{Y}/G$ is ramified 
where any of the roots, or their product, is equal to zero.  
There are $(\dim \mathcal{G} - r)$ roots, so this defines a 
hypersurface of degree $(\dim \mathcal{G} - r)$ in the total space of
$\omega_{Y}$.  The ramification divisor of $X \to Y$ is the 
intersection of this hypersurface with the section $\phi$,
which is the divisor 
corresponding to the line bundle $\omega_{Y}^{\otimes(\dim \mathcal{G} - r)}.$
Thus the degree of the branch divisor will be $(\dim \mathcal{G} - r)(2g-2)$.

Combining all of this information, we see that the dimension of the 
Prym is
\begin{eqnarray*}
\dim \Prym_{\mathbb{t}} (X) 
& = & r (g - 1) + (r-(r-1)) \frac{(\dim \mathcal{G} - r)(2g-2)}{2} \\
& = & r(g-1)  + (\dim \mathcal{G} - r)(g-1) \\
& = & \dim \mathcal{G}(g-1).
\end{eqnarray*}

By comparison, the dimension of the base space is
\begin{equation*}
\Sigma_{i=1}^{r} h^{0}(Y,\omega_{Y}^{d_i})
\end{equation*}

The sum of the degrees $d_{i}$ of the basic invariant polynomials of 
$\mathfrak{g}$ is the dimension of a Borel subalgebra, $(\dim 
\mathcal{G} + r)/2$.  For $g>1$, Riemann-Roch gives 
\begin{eqnarray*}
\Sigma_{i=1}^{r} h^{0}(Y,\omega_{Y}^{d_i})
& = & \Sigma_{i=1}^{r} (2d_{i}-1) (g-1)  \\
& = & (\dim \mathcal{G} + r - r)(g-1) \\
& = & \dim \mathcal{G}(g-1).
\end{eqnarray*}

Which, as Hitchin said, ``somewhat miraculously'' turns out to be the same thing.  

Markman [Ma] and Bottacin [B] generalized the Hitchin system by twisting
the line bundle $\omega_{Y}$ by an effective divisor $D$.  The effect of
this is to create a family of integrable systems, parametrized by 
the residue of the Higgs field $\phi$ at $D$.  The base space of each 
system is a fiber of the map 
\begin{gather*}
B := \bigoplus_{i=1}^{r}H^{0}(Y,\omega_{Y}(D)^{\otimes d_i}) \\
\downarrow \\ 
\bar{B} :=  \text{the space of possible residues at }  D
\end{gather*}
which sends the set of $r$ sections in $B$ 
to its set of residues at $D$.  At each point of $D$, there are $r$ independent 
coefficients, so the dimension of $\bar B$ is $r(\deg D)$.  Thus the base space 
of each system has dimension
\begin{eqnarray*}
\dim B - \dim \bar B
& = & \sum_{i=1}^{r}h^{0}(Y,\omega_{Y}(D)^{\otimes d_i}) - r(\deg D) \\ 
& = & \sum_{i=1}^{r}( d_{i}(2g-2 + \deg D) - (g-1)) -r(\deg D) \\
& = & (1/2)(\dim \mathcal{G} + r)(2g - 2 + \deg D) -r(g-1) -r(\deg D) \\
& = & (\dim \mathcal{G})(g-1) + \frac{\dim \mathcal{G} -r}{2} \deg D  
\end{eqnarray*}

Markman showed that the generic fiber of this system is again isogenous to 
$\Prym_{\mathfrak{t}}(X)$, where $X$ is a cameral cover of the base curve $Y$.  
The construction of the cameral cover is similar to the case of the 
Hitchin system, except that $\phi$ is a section of $ad(P) \otimes 
\omega_{Y}(D)$.  Thus the ramification divisor is 
$(\omega_{Y}(D))^{\otimes(\dim \mathcal{G} -r)}$,  and the dimension is
\begin{eqnarray*}
\dim \Prym_{\mathfrak{t}} (X) 
& = & r (g - 1) + \frac{(\dim \mathcal{G} - r)(2g-2 + 
\deg D)}{2} \\
& = & \dim \mathcal{G}(g-1) + \frac{(\dim \mathcal{G} -r)}{2} \deg D.
\end{eqnarray*}

Again, this is the same dimension as the base of the system.


\begin{thebibliography}{BBB}
   \bibitem[A]{A}
	P. Aspinwall, \emph{Aspects of the Hypermultiplet Moduli Space
	in String Duality}, hep-th/9802194.
   \bibitem[AvM]{AvM}
	M. Adler and P. van Moerbeke, \emph{Completely integrable systems,
	Euclidean Lie algebras, and curves}, Adv. in Math. \textbf{38} 
	(1980), no. 3, 267-317.
   \bibitem[B]{B}
        F. Bottacin, \emph{Symplectic geometry on moduli spaces of stable 
	pairs}, Ann. Sci. ƒcole Norm. Sup. (4) \textbf{28} (1995), no. 4,
	391--433.
   \bibitem[BZ]{BZ}
	Ya. G. Berkovich and E. M. Zhmud', \emph{Characters of Finite Groups.
	Part 1}, translated by P Shumyatsky and V. Zobina, translation edited
	by David Louvish.  American Mathematical Society, Providence, Rhode
	Island, 1998.
   \bibitem[D1]{D1}
	R. Donagi, \emph{Seiberg-Witten Integrable Systems}, Algebraic 
        geometry---Santa Cruz 1995, 3--43, Proc. Sympos. Pure Math., 
	\textbf{62}, Part 2, Amer. Math. Soc., Providence, RI, 1997.  
	alg-geom/9705010.
   \bibitem[D2]{D2}
	R. Donagi, \emph{Decomposition of spectral covers}, in Journees de
        Geometrie Algebrique D'Orsay, Asterisque \textbf{218} (1993), 145-175.
   \bibitem[F]{F}
	G. Faltings, \emph{Stable G-bundles and Projective Connections}, 
	Journal of Algebraic Geometry \textbf{2} (1993), 507-568.
   \bibitem[FH]{FH}
	W. Fulton and J. Harris, \emph{Representation Theory:  A First
	Course}, Springer-Verlag, New York, 1991.
   \bibitem[H]{H}
	N. J. Hitchin, \emph{Stable bundles and integrable systems}, Duke
	Math. J. \textbf{54}, no. 1, 91-114.
   \bibitem[Ma]{Ma}
        E. Markman, \emph{Spectral curves and integrable systems}, Comp. 
        Math. \textbf{93} (1994), 255-290.
   \bibitem[Me]{M}
	J.-Y. M\'{e}rindol, \emph{Vari\'{e}t\'{e}s de Prym d'un 
	rev\^{e}tement gaoloisien}.  J. reine angew. Math. \textbf{461} (1995), 
	49-61.
   \bibitem[MS]{MS}
	A. McDaniel and L. Smolinsky, \emph{A Lie Theoretic Galois Theory
	for the Spectral Curves of an Integrable System. II}.  Transactions
	of the American Mathematical Society \textbf{49}, no. 2 (February 1997),
	747-762.  
   \bibitem[S]{S}
	R. Scognamillo, \emph{An elementary approach to the abelianization
	of the Hitchin system for arbitrary reductive groups}.  Comp.
	Math. \textbf{110} (1998), 17-37.
\end{thebibliography}
\end{document}